\documentclass[twoside, 10pt]{amsart}

\def\docTitle{The resolution property of algebraic surfaces}
\def\docAuthor{Philipp Gross}
\def\docSubject{Algebraic Geometry}
\def\docKeywords{enough locally frees, resolution property, global resolutions, vector bundles, singular surfaces, non-projective surfaces, gluing schemes, divisorial schemes, AF-schemes, quasiprojective open subsets}

\title{The resolution property of algebraic surfaces}
\author[P. Gross]{Philipp Gross}
\address[P. Gross]{Mathematisches Institut\\ Heinrich-Heine-Universit\"at\\ D-40225 D\"usseldorf\\ Deutschland}
\email{gross@math.uni-duesseldorf.de}
\thanks{The author was supported by the Deutsche Forschungsgemeinschaft, Forschergruppe 790 ``Classification of Algebraic Surfaces and Complex Manifolds''.}
\date{\today}

\usepackage{ifpdf}
\ifpdf
\usepackage[pdfstartview = {FitH},
  plainpages={false},
  pdfpagelabels={true}, 
  pdfpagemode={UseOutlines},
  bookmarksopen = {false},
  colorlinks = {false}, linkcolor = {blue}, citecolor = {green},
  pdftitle={\docTitle}, pdfauthor={\docAuthor}, pdfsubject={\docSubject}, pdfkeywords={\docKeywords} ]{hyperref}%
\fi

\usepackage{srcltx}  				
\usepackage{amsthm}
\usepackage{amsmath}
\usepackage{amssymb,amsfonts}	
\usepackage{mathrsfs}
\usepackage[mathscr]{eucal}
\usepackage{fixltx2e}[2005/12/01]
\usepackage[all]{xy}	

\usepackage{enumitem}
\setlist{itemsep=0.5em}
\usepackage{natbib}				
\bibpunct{[}{]}{,}{n}{}{;}			

\newtheoremstyle{dtheorem}{3 mm}{1 mm}{\itshape}{}{\bfseries}{.}{ }
  {\thmnumber{(#2) }\thmname{#1}\thmnote{ \mdseries(#3)\bfseries}}
\newtheoremstyle{ddef}{3 mm}{1 mm}{\normalfont}{}{\bfseries}{.}{ }
  {\thmnumber{(#2) }\thmname{#1}\thmnote{ \mdseries(#3)\bfseries}}
\newtheoremstyle{dremark}{3 mm}{1 mm}{\normalfont}{}{\itshape}{}{ }
  {\thmnumber{\upshape\bfseries(#2) }\itshape\mdseries\thmname{#1}.\thmnote{\;\mdseries #3 ---}}

\theoremstyle{dtheorem}
\newtheorem{thm}{Theorem}[section]

\numberwithin{equation}{thm}
\makeatletter
\@addtoreset{thm}{section}
\makeatother

\newtheorem*{thm*}{Theorem}
\newtheorem{prop}[thm]{Proposition}
\newtheorem{lem}[thm]{Lemma}
\newtheorem{cor}[thm]{Corollary}

\theoremstyle{ddef}
\newtheorem{ex}[thm]{Example}
\newtheorem{defi}[thm]{Definition}

\newtheorem*{claim*}{Claim}
\theoremstyle{dremark}
\newtheorem{rem}[thm]{Remark}

\newtheorem{problem}[thm]{Problem}

\newcommand{\msheaf}[1]{\mathcal{#1}}		

\newcommand{\shC}{\msheaf{C}}

\newcommand{\shE}{\msheaf{E}}
\newcommand{\shEnd}{\ensuremath{ \msheaf{E}nd }}
\newcommand{\shExt}{\msheaf{E}xt}
\newcommand{\shF}{\msheaf{F}}
\newcommand{\shG}{\msheaf{G}}
\newcommand{\shH}{\msheaf{H}}
\newcommand{\shHom}{\msheaf{H}om}
\newcommand{\shI}{\msheaf{I}}

\newcommand{\shL}{\msheaf{L}}
\newcommand{\shM}{\msheaf{M}}
\newcommand{\shN}{\msheaf{N}}

\newcommand{\shO}{\msheaf{O}}

\newcommand{\shS}{\msheaf{S}}

\newcommand{\bbC}{\mathbb{C}}

\newcommand{\bbN}{\mathbb{N}}

\newcommand{\bbP}{\mathbb{P}}
\newcommand{\bbQ}{\mathbb{Q}}

\newcommand{\bbZ}{\mathbb{Z}}

\newcommand{\ideal}[1]{\mathfrak{#1}}

\newcommand{\coloneq}{\ensuremath{:=}}

\newcommand{\AF}{\text{AF}}

\newcommand{\red}{\text{red}}

\DeclareMathOperator{\coker}{\operatorname{coker}}

\DeclareMathOperator{\im}{\operatorname{im}}

\DeclareMathOperator{\codim}{\operatorname{codim}}
\DeclareMathOperator{\Coho}{\operatorname{H}}

\DeclareMathOperator{\Ext}{\operatorname{Ext}}
\DeclareMathOperator{\Hom}{\operatorname{Hom}}
\DeclareMathOperator{\pd}{\operatorname{pd}}

\DeclareMathOperator{\Ann}{\operatorname{Ann}}		
\DeclareMathOperator{\Ass}{\operatorname{Ass}}

\DeclareMathOperator{\trdeg}{\operatorname{deg.tr}}
\DeclareMathOperator{\height}{\operatorname{ht}}

\DeclareMathOperator{\Spec}{\operatorname{Spec}}
\DeclareMathOperator{\Supp}{\operatorname{Supp}}

\DeclareMathOperator{\Pic}{\operatorname{Pic}}

\newcommand{\pf}[1]{ {#1}_* }
\newcommand{\pb}[1]{ {#1}^* }

\newcommand{\us}[1]{ {#1}^!}

\newcommand{\ses}[3]{\ensuremath{ 
  0 \rightarrow #1 \rightarrow #2 \rightarrow #3 \rightarrow 0
}}

\newcommand{\xses}[5]{\ensuremath{ 
  0 \rightarrow #1 \overset{#2}\rightarrow #3 \overset{#4}\rightarrow #5 \rightarrow 0
}}

\newcommand{\es}[3]{\ensuremath{ 
  #1 \rightarrow #2 \rightarrow #3 
}}\newcommand{\xes}[5]{\ensuremath{ 
  #1 \overset{#2}\rightarrow #3  \overset{#4}\rightarrow #5 
}}

\begin{document}

\begin{abstract}
We prove that on separated algebraic surfaces any coherent sheaf is a quotient of a locally free sheaf. This class contains many schemes that are neither normal, reduced, quasiprojective or embeddable into toric varieties.
Our methods extend to arbitrary $2$-dimensional schemes that are proper over a noetherian ring.
\end{abstract}

\maketitle

\tableofcontents

\section*{Introduction} 

A noetherian scheme  (or complex analytic space, or more generally a locally ringed site) has the \emph{resolution property} or \emph{enough locally free sheaves} if every coherent sheaf $\shM$ admits a surjection $\shE \twoheadrightarrow \shM$ by some coherent locally free sheaf $\shE$ (also called a vector bundle). 
For an introduction to this property we refer to the seminal paper of Totaro \cite{MR2108211}.

The aim of this article is to prove the resolution property for all separated algebraic surfaces, generalizing a result of Schr\"oer and Vezzosi who verified the resolution property for \emph{normal} separated algebraic surfaces  \cite{MR2041778} (see Theorem \ref{THM:main_theorem} and Corollary \ref{COR:algebraic_surfaces}).
Our methods extend to arbitrary $2$-dimensional schemes that are proper over an arbitrary noetherian base ring, throughout denoted by $A$. 

\medskip

Unless stated otherwise all schemes are assumed to be noetherian. The word \emph{surface} refers to a $2$-dimensional separated $A$-scheme  of finite type. 

\medskip

It is known that a scheme satisfies the resolution property if it has an ample line bundle, or more generally, if it has an ample family of line bundles; that is a family of invertible sheaves where the whole collection behaves like an ample line bundle (see \cite[Expos\'e II, 2.2]{sga6}  and \cite{MR1970862}). 
This includes all schemes that are quasiprojective over a noetherian ring and all $\bbQ$-factorial schemes with affine diagonal  \cite{MR1970862}. The case of regular, separated schemes is also known as Kleiman's Theorem \cite{MR0219545} (independently proven by Illusie \cite{sga6}).

There exist regular algebraic schemes that do not have the resolution property \cite{MR2041778}. 
However, they are not separated because they do not have affine diagonal. The latter means that the intersection of two open affines is affine which is a necessary condition for the resolution property to hold by Totaro \cite[1.3]{MR2108211}.

The situation is different in the analytic category. Schuster showed that every smooth, compact, complex surface satisfies the resolution property \cite{MR676049} whereas
it fails for generic complex tori of dimension $\geq 3$ \cite[A.5]{MR1902630}.

\medskip

Starting in dimension $2$, it can happen that a scheme has no effective Cartier-Divisors or even no non-trivial line bundles at all \cite{MR1726231}. Therefore it is not possible to construct locally free resolutions by invertible sheaves, in general. The lack of control for vector bundles of higher rank makes it difficult to tackle the resolution property for singular and non-projective schemes.
To the authors knowledge, it is not known if there exists a scheme with affine diagonal that does not have the resolution property. 
Even the case of normal, toric, separated, algebraic schemes of dimension $\geq 3$ is completely open \cite{MR2448277}. 
By our work we seek to give a clarifying picture of the resolution property for separated schemes of dimension $\leq 2$.

\medskip

The paper is organized as follows. In the first section we shall investigate quasiprojectivity of large open subsets. Given a separated $A$-scheme of finite type we prove that every point has a quasiprojective neighborhood whose complement has codimension $\geq 2$ (Theorem \ref{THM:thick_quasiprojectives}).
This is well-known if $X$ is normal \cite{MR2146196}. Indeed, since $X$ admits only trivial finite birational covers in that case, this is a consequence of Chow's Lemma and Zariski's Main Theorem.
We will deduce the general case using Ferrand's pinching techniques \cite{MR2044495} and deformation theory of vector bundles \cite[5.A]{MR2223409}.
Consequently, every coherent sheaf is a quotient of a coherent sheaf which is locally free outside a closed subset of codimension $\geq 2$. 
Restricting to surfaces we generalize the gluing methods of Schr\"oer and Vezzosi in Section \ref{SEC:glueing_resolutions}
to construct locally free resolutions of coherent sheaves that are already locally free outside finitely many points of codimension $2$. The cohomological obstructions appearing here lie in second cohomology groups of coherent sheaves.
In order to control these, we study in  Section \ref{SEC:cample} a partial cohomological vanishing condition for families of coherent sheaves $\shE_n$, $n \in \bbN$, on a proper scheme $X$ of any dimension. Given an integer $1 \leq d < \dim(X)$ we call $(\shE_n)$ $d$-ample if for every coherent sheaf $\shM$ the groups $\Coho^i(X,\shM \otimes \shE_n)$ vanish for all $i > d$. Here our main result states that $(\dim(X)-1)$-ampleness is preserved and reflected under pullback by alterations.
We use this in Section \ref{SEC:existence_of_cample_families_of_vectorbundles} to construct a $1$-ample family of vector bundles of rank $2$ on an arbitrary proper surface.
In the last section we collect the preceding results to prove the resolution property for a large class of $2$-dimensional separated schemes.

\medskip

\subsection*{Acknowledgments}

This article consists of a part of my Ph.D. thesis. I am indebted
to my advisor Stefan Schr\"oer for his encouragement to pursue this project and for teaching me algebraic geometry. 
I would also like to thank Jarod Alper, Christian Liedtke, Holger Partsch, David Rydh, Sasa Novakovic and Felix Sch\"uller for many inspiring conversations and helpful suggestions. 

Finally, I thank the Deutsche Forschungsgemeinschaft for their financial support as a member of the Forschergruppe 790 ``Classification of Algebraic Surfaces and Complex Manifolds''.

\section{Thick quasiprojective open subsets}
\label{SEC:thick_quasiprojectives}

The aim of this section is to provide an existence result of thick quasi-projective open subschemes of a scheme $X$ that is separated and of finite type over $A$ (see Theorem \ref{THM:thick_quasiprojectives}).
Here we call a subset $V \subset X$ \emph{thick} if it is open in $X$ and $\codim(X-V,X) \geq 2$. Clearly, the finite intersection of thick subsets is thick. If $W \subset V$, $V \subset X$ are thick then $W$ is thick in $X$.

We will frequently use birational auxiliary schemes $Y$ that are quasiprojective over $A$ and are endowed with a birational map $Y \to X$. Let us call a morphism of schemes $f \colon Y \to X$  ($U$-\emph{admissible}) \emph{birational} if there exists dense open subschemes $U \subset X$, $f^{-1}(U) \subset Y$  such that the restriction $f_U$ is an isomorphism.

We start with a sequence of preparatory lemmas.
\begin{lem}\label{LEM:birat_quasifinite_in_codim1}
 Let $f \colon Y \to X$ be a birational morphism of finite type of locally noetherian schemes. Then $f$ is quasifinite over all points of codimension $\leq 1$.
\end{lem}
\begin{proof}
Let $x \in X$ be a point of codimension $1$. By applying base change with $\Spec \shO_{X,x} \to X$ we may assume that $X$ is local of Krull dimension $1$ with closed point $x$. Moreover, we may assume that $X$ and $Y$ are irreducible by a second base change.
Let $y \in f^{-1}(x)$ a maximal point. Then $\dim \shO_{Y,y} \geq 1$ because $y$ cannot be a generic point of $Y$, whereas $0 \leq \dim \shO_{Y,y} + \trdeg_{k(x)}k(y) \leq \dim \shO_{X,x} = 1$  by \cite[5.6.5.1]{egaIV_2}. Consequently, holds $\trdeg_{k(x)}k(y)=0$ which shows that $f^{-1}(x)$ is finite and discrete.
\end{proof}

\begin{lem}\label{LEM:codim_qf_maps}
Let $f \colon Y \to X$ be a quasifinite map of finite type of locally noetherian schemes. 
Then $\codim(f(Z),X) \geq \codim(Z,Y)$ for every closed subset $Z \subset Y$.
\end{lem}
\begin{proof}
For each $z \in Z$ holds \(\dim \shO_{f(z)}  \geq  \dim \shO_{z} - \dim\shO_{z} \otimes_{\shO_{f(x)}} k(f(x)) = \dim \shO_{z} \)
by  \cite[5.4.2]{egaIV_2} so that
\[ \codim(f(Z),X) = \inf_{z \in Z} \dim \shO_{f(z)}  \geq \inf_{z \in Z} \dim \shO_{z} = \codim(Z,X). \qedhere\]
\end{proof}

\begin{lem}\label{LEM:properbirat_finite_over_thick}
Let $f \colon Y \to X$ be a proper birational morphism of locally noetherian schemes which is $U$-admissible for a dense open subset $U \subset X$. 
Then $f$ is finite over a thick  subset $V \subset X$ containing $U$. 
\end{lem}
\begin{proof}
Let $Y \overset{f'} \to X' \overset{g}\to X$ be a Stein factorization of $f$. 
From Zariski's Main Theorem \cite[4.4.1]{egaIII_1} and  Lemma \ref{LEM:birat_quasifinite_in_codim1} we infer that  $f'$ is an isomorphism over a thick open subset $V' \subset X'$. Then $Z \coloneq g(X'-V') \subset X$ has codimension $\geq 2$ by Lemma \ref{LEM:codim_qf_maps} since $g$ is finite. Since $g^{-1}(X-Z) = X' - g^{-1}(Z) \subset X' - (X'-V') = V'$ we conclude that $f$ is finite over $X-Z$.
\end{proof}

\begin{lem}\label{LEM:properbirat_iso_over_thick}
Let $f \colon Y \to X$ be a proper birational morphism of locally noetherian schemes which is $U$-admissible for a dense open subset $U \subset X$. Suppose that $\shO_{X,x}$ is normal (equivalently regular) for all $x \in X-U$ with $\dim \shO_{X,x}=1$. 
Then $f$ is $V$-admissible for a thick  subset $V \subset X$ containing $U$.
\end{lem}
\begin{proof}
By Lemma \ref{LEM:properbirat_finite_over_thick} we may assume that $f$ is finite by replacing $X$ with a suitable thick subset. Then for all $x \in X$ with $\dim \shO_{X,x} \leq 1$ the map $\shO_{X,x} \to (\pf f \shO_Y)_x$ is an isomorphism. This is obvious if $x \in U$ but if $x \in X-U$ then $\shO_{X,x}$ is normal by hypothesis and the assertion follows from \cite[4.4.9]{egaIII_1}. Thus $\shO_X \to \pf f \shO_Y$ is an isomorphism outside  a closed subset $Z\subset X$ of codimension $\geq 2$ and it follows $f$ is an isomorphism over $X-Z$.
\end{proof}

Recall that a scheme is called an $\AF$-scheme if every finite set of points is contained in an affine open neighborhood. This is also known as the Kleiman-Chevalley condition \cite{MR0206009} and beyond that related to embeddings into toric varieties \cite{MR1227474} or to \'etale cohomology \cite[\S 4]{MR0289501}.
It is always satisfied if there exists an ample line bundle \cite[4.5.4]{egaII}; for example, if $X$ is quasiprojective over a noetherian ring. The converse holds for proper, smooth algebraic schemes (see \cite{MR0206009} and its generalization by W{l}odarczyk \cite{MR1683254})  but not for non-normal schemes.
\begin{rem}\label{REM:Ferrand_AF} 
Suppose we have a finite, birational map with schematically dense image of noetherian schemes $f \colon X' \to X$. The latter means $\shO_{X} \to \pf f \shO_{X'}$ is injective which is automatic if $X$ is reduced.
Let $Y \subset X$ be the \emph{conductor subscheme}, defined by the conductor ideal $\Ann \coker (\shO_X \to \pf f \shO_{X'}) \subset \shO_X$ and define $g \coloneq f_Y \colon f^{-1}(Y) \to Y$. Then $X$ is isomorphic to the pinching $Y \coprod_g X'$ of $Y$ with $X'$ along $g$ by Ferrand \cite[4.3]{MR2044495}. 
By Ferrand's Theorem \cite[5.4]{MR2044495} follows that $X$ is an $\AF$-scheme if and only if $X'$ and $Y$ are $\AF$-schemes. 
The latter is satisfied if $X'$ and $Y'$ are quasiprojective over a noetherian ring but then it does not follow that $X$ is quasiprojective (see Example \ref{EX:nondivisorial_surface} below).
\end{rem}

\begin{thm}\label{THM:thick_quasiprojectives}
Let $X$ be a separated scheme of finite type over a noetherian ring $A$. Then every point $x \in X$ has a thick neighborhood $x \in V \subset X$ that is quasiprojective over $A$.
\end{thm}
\begin{proof}
Assume first that $X$ is reduced.
Let $x \in U \subset X$ be an affine open set. By enlarging $U$ with disjoint affine open sets we may assume that $U$ is dense. Then by Nagata there exists a $U$-admissible blow up $f \colon X' \to X$ such that $X'$ is quasiprojective over $A$ \cite[2.6]{MR2356346}. 

By Lemma \ref{LEM:properbirat_finite_over_thick} we may assume that $f$ is finite by replacing $X$ with a thick subset. Moreover, $f$ is birational and has schematically dense image.
Let $Y \subset X$ be conductor subscheme of $f$.
Choose a closed subset $Z \subset Y$ such that $Y-Z \subset Y$ is a dense, affine and hence a quasiprojective open subscheme. 
Then $\codim(Z,X) \geq \codim(Z,Y) + \codim(Y,X) \geq 2$. So by replacing $X$ with $X-Z$ we may assume that $Y$ and $X'$ are quasiprojective over $A$. 

In light of Remark  \ref{REM:Ferrand_AF} we may therefore assume that $X$ is an $\AF$-scheme.
In particular, there exists a dense affine open neighborhood $x \in U_1 \subset X$ that contains the finite set of points $z \in X$ of codimension $1$ with $\shO_{X,z}$ non-regular. So by repeating the previous arguments with $U$ replaced by $U_1$, we may assume that for all $z \in X-U$ with $\dim \shO_{X,z} \leq 1$ the stalk $\shO_{X,z}$ is regular. 
But then $f$ is an isomorphism over a thick  subset by Lemma \ref{LEM:properbirat_iso_over_thick}.

Let $X$ now be arbitrary. Replacing $X$ with a thick subset we may assume that $X_\red$ is quasiprojective over $A$ by the special case.
By removing closed subsets $Z \subset X - U$ with $\codim(Z,X-U) \geq 1$ we may assume that $X -U$ is affine so that $X = U \cup (X-U)$ has cohomological dimension $\leq 1$ \cite[2.8]{MR2096147}. Therefore every ample $\shO_{X_\red}$-module lifts to an invertible $\shO_X$-module $\shL$ since the obstructions lie in $2^{nd}$ cohomology groups of coherent sheaves \cite[Theorem 5.3]{MR2223409}. Now $\shL$ is also ample by \cite[4.5.13]{egaII} and we conclude that $X$ is quasiprojective over $A$.
\end{proof}

We derive an existence result for affine open neighborhoods, generalizing a result of Raynaud \cite[VIII 1]{MR0260758}.
\begin{cor}\label{COR:points_of_codim_1_are_contained_in_affine_open}
Let $X$ be a separated scheme of finite type over a noetherian ring. Then every finite set of points $x_1, \dots, x_n \in X$ is contained in a dense affine open neighborhood as long as $x_2, \dots, x_n$ are of codimension $\leq 1$, .
\end{cor}
\begin{proof}
By Theorem \ref{THM:thick_quasiprojectives} there exists a neighborhood $x_1 \in V \subset X$ that is thick and quasiprojective. But then also $x_i \in V$ for $i \geq 2$ since $\codim(X-V,X) \geq 2$ and the result follows using that $V$ is an $\AF$-scheme.
\end{proof}
The bound on the codimension in Theorem \ref{THM:thick_quasiprojectives} and Corollary \ref{COR:points_of_codim_1_are_contained_in_affine_open} is sharp because there exist irreducible algebraic surfaces having a pair of closed points that do not admit an affine open neighborhood \cite{MR1726231}.

\medskip

In the proof of Theorem \ref{THM:thick_quasiprojectives} we have seen that up to a closed subset of codimension $\geq 2$ the scheme $X$ is the pinching of quasiprojective schemes.
The following example illustrates that pinching may destroy many Cartier divisors, so that $X$ is a priori not even \emph{divisorial}; that is, the complements of effective Cartier divisors define a base for the Zariski topology (\cite[Expos\'e II, 2.2]{sga6},\cite{MR0219545}).

\begin{ex}[A non-divisorial proper algebraic surface whose normalization is projective]\label{EX:nondivisorial_surface}
We work over an algebraically closed field $k$, say $k \coloneq \bbC$ for simplicity. 
Let $E$ be an elliptic curve and consider the surface $X \coloneq E \times \bbP^1$ with projections $p \colon X \to E$ and $q \colon X \to \bbP^1$. Choose distinct fibers $E_0$ and $E_\infty$ over $\bbP^1$.
Let $t_x \colon E \to E$ be the translation with respect to a rational point $x\in E$ of infinite order. 
Define the finite map $g \colon E_0 \coprod E_\infty \to E$ as the identity on $E_0$ and as $t_x$ on $E_\infty$. 

By Ferrand \cite[5.4]{MR2044495}, the pushout $S$ of the closed immersion $ i \colon E_0 \coprod E_\infty \hookrightarrow X$ along $g$ exists in the category of schemes and fits in cartesian and cocartesian square
\[
\xymatrix{ E_0 \coprod E_\infty \ar[r]^i\ar[d]^g & X \ar[d]^f \\
 E \ar[r]^j & S}
\]
Here $j$ is a closed immersion and $f$ is finite with schematically dense image. It follows that $S$ is an integral, proper surface with normalization $f$.

Let us see that $S$ is not divisorial by way of contradiction. Assume that for a given point $y \in E \subset S$ there exists an effective Cartier divisor $C \subset S$ with $y \notin C$ and $S-C$ affine. Then $C \cap E$ is non-empty and $0$-dimensional. It follows that the line bundle $\shL \coloneq \shO_S(C)|_E$ has positive degree and hence is ample.  Now the natural isomorphism $\pb g \shL \simeq \pb i \pb f \shO_S(C)$ induces the isomorphisms
$\shL \simeq \pb g \shL|_{E_0} \simeq \pb f \shO_S(C)|_{E_0}$ and $\pb f \shO_S(C)|_{E_\infty}\simeq \pb {t_x} \shL$.

Since $X$ is a ruled surface holds $\pb f \shO_S(C) \simeq \pb p \shM \otimes \pb q \shO_{\bbP^1}(n)$ for some $\shM \in \Pic(E)$ and $n \in \bbN$. It follows that $\pb f \shO_S(C)|_{E_\infty} \simeq \shM \simeq \pb f \shO_S(C)|_{E_0}$ and consequently $\shL \simeq \pb{t_x} \shL$.
 But then $x$ must have finite order by the theory of abelian varieties \cite[p. 60, Application 1]{MR0282985}, contradicting the choice of $x$.
\end{ex}

As a consequence of the existence of thick quasiprojective open neighborhoods we obtain resolutions by coherent sheaves that are invertible outside a closed subset of codimension $\geq 2$ by extending coherent sheaves beyond dense affine open subsets.
\begin{prop}\label{PROP:almost_ample_resolver}
Let $X$ be a separated scheme that is of finite type over $A$. Then for every $x \in X$ there exists a coherent sheaf $\shF$ with the following properties:
\begin{enumerate}
\item\label{PROP:almost_ample_resolver:IT:1} For every coherent sheaf $\shM$ and every $m \gg0$ there exists a map \[(\shF^{\otimes m})^{\oplus n} \rightarrow \shM\] for some $n \in \bbN$ which is surjective near $x \in X$.
\item \label{PROP:almost_ample_resolver:IT:2} There exists a thick quasiprojective neighborhood $x \in V \subset X$ such that 
$\shF|_V$ is invertible and $\shF|_V^{\vee}$ is an ample $\shO_V$-module.
\item \label{PROP:almost_ample_resolver:IT:3} There exists a $V$-admissible blow-up $f \colon X' \to X$ such that $\pb {f_V} \shF|_V^{\vee}$ extends to an ample $\shO_{X'}$-module.
\end{enumerate} 
\end{prop}
\begin{proof}
By Theorem \ref{THM:thick_quasiprojectives} there exists a thick neighborhood $x \in V \subset X$ that is quasiprojective over $A$.
Let $f \colon X' \to X$ be a $V$-admissible blow-up such that $X'$ is quasiprojective over $A$.
Choose an ample $\shO_V$-module $\shL$ such that $\pb {f_V} \shL$ extends to an ample $\shO_{X'}$-module.
Let $x \in U \subset V$ be an affine open neighborhood with $U=V_s$ for some $s \in \Coho^0(V, \shL^n)$ and $n > 0$.
Let $\shI \subset \shO_X$ be a coherent ideal with $V(\shI)=X-U$. Then the isomorphism of $\shO_U$-modules $\shO_{U} \xrightarrow{\sim} \shI|_{U}$ extends to a map of $\shO_V$-modules $\shO_V \to \shI|_V \otimes \shL^{nm}$ for some $m >0$ \cite[6.8.1]{egaI2nd}. Hence the twist $\shL^{-nm} \to \shI|_V$ is an isomorphism near $x$. 
Then $\shL^{-mn} \to \shI|_V$ extends to a map $\varphi \colon \shF \to \shI$ for some coherent $\shO_X$-module $\shF$ with $\shF|_V \simeq \shL^{-mn}$.
It follows that \ref{PROP:almost_ample_resolver:IT:2} and \ref{PROP:almost_ample_resolver:IT:3}  are satisfied. 

In order to prove \ref{PROP:almost_ample_resolver:IT:1} let $\shM$ be  an arbitrary coherent sheaf on $X$. 
Then every surjection $\shO_{U}^{\oplus n} \twoheadrightarrow \shM|_{U}$ extends to a map $(\shI^m)^{\oplus n} \to \shM$ for all  $m \gg 0$ \cite[I.6.9.17]{egaI2nd}.
Using that the maps $\varphi^{\otimes m} \colon \shF^{\otimes m} \to  \shI^{\otimes m}$ and $\shI^{\otimes m} \twoheadrightarrow \shI^m$ are surjective near $x$, we infer that the same holds for the composition 
\[(\shF^{\otimes m})^{\oplus n} \to (\shI^{\otimes m})^{\oplus n} \to (\shI^{m})^{\oplus n} \to \shM. \qedhere\]
\end{proof}

\begin{defi}
We call coherent sheaf $\shF$ \emph{almost anti-ample (near $x \in X$)} if it satisfies properties \ref{PROP:almost_ample_resolver}.\ref{PROP:almost_ample_resolver:IT:1}-\ref{PROP:almost_ample_resolver:IT:3}.
\end{defi}

In case that $X$ has dimension $2$, an almost anti-ample sheaf is invertible outside finitely many points of codimension $2$.

\section{Gluing resolutions}
\label{SEC:glueing_resolutions}

In this section we formulate conditions that are sufficient for the existence of locally free resolutions of coherent sheaves which are locally free away from finitely many closed points of codimension $2$. 

We pursue the strategy of constructing surjections $\varphi \colon \shE \twoheadrightarrow \shM$ with predefined kernel $\shS$, the  \emph{first syzygy} of $\varphi$, generalizing the methods of Schr\"oer and Vezzosi \cite{MR2041778}. Such a map $\varphi$ is determined by an extension class $\gamma \in \Ext^1(\shM, \shS)$ up to isomorphism of $\shE$ over $\shM$. Instead of gluing morphisms one glues extension classes. 
This is controlled by the exact sequence
\begin{equation*}
\xes{\Ext^1(\shM,\shS)}{\ell}{\Coho^0(X, \shExt(\shM,\shS))}{}{ \Coho^2(X,\shHom(\shM,\shS))}
\end{equation*}
which can be read off the local to global spectral sequence for $\Ext$. 
Here the map $\ell$ maps a class $[e]$ represented by an extension $e \colon \ses{\shS}{\shE}{\shM}$ to a global section $\ell(e)$ whose stalk $\ell(e)_x$ at $x \in X$ is given by the class of the localized sequence $[e_x] \in \Ext^1_{\shO_{X,x}}(\shM_x, \shS_x) \simeq \shExt^1_{\shO_{X,x}}(\shM_x,\shS_x)_x$.

\medskip

If $\shM$ is locally free outside finitely many closed points $Z \subset X$, the local extensions of $\shM$ by $\shS$ appear in a simple form. In fact, there is a canonical isomorphism
$ \bigoplus_{z \in Z} \Ext^1_{\shO_{X,z}}(\shM_z, \shS_z) \simeq \Coho^{0}(X,\shExt^1_{\shO_X}(\shM,\shS))$. This enables us to choose elements of $\Ext^1_{\shO_{X,z}}(\shM_z, \shS_z)$, $z \in Z$, independently to a get a global section of $\shExt^1_{\shO_X}(\shM,\shS)$.

\medskip

We will use a more convenient formulation of the obstruction space for gluing local resolutions of $\shM$ by $\shS$. 
Note that $\shHom(\shM, \shO_X)$ and $\shHom(\shM, \shS)$ are isomorphic over the open subset where $\shM$ is locally free. In general, they are not isomorphic everywhere since the tensor product may have torsion sections.
Nevertheless, if $\shM$ is locally free up to finitely many closed points, we deduce from the succeeding Lemma \ref{LEM:coho_almost_isom} the isomorphism on the top cohomology
\begin{equation}
\Coho^2(X, \shHom(\shM, \shS)) \simeq \Coho^2(X, \shM^\vee \otimes \shS).
\end{equation}
Analogously, the obstruction space does not change up to isomorphism if $\shS$ is modified over finitely many closed points.

\begin{lem}\label{LEM:coho_almost_isom}
Let $X$ be a noetherian scheme and $\shF$, $\shF'$ be coherent sheaves that are isomorphic outside a closed subscheme $Z \subset X$. Then for all $i \in \bbN$ with $i \geq \dim(Z)$ holds \(\Coho^{i+2}(X, \shF) \simeq \Coho^{i+2}(X, \shF') \) .
\end{lem}
\begin{proof}
If there exists a map $ u \colon \shF \to \shF'$ that is an isomorphism on $X-Z$ then the support of  $\ker u$ and $\coker u$ has dimension $\leq d$. So the assertion follows by taking cohomology of the two exact sequences
\[
 \ses{\ker u}{\shF}{\im u} \quad \text{and} \quad \ses{\im u}{\shF'}{\coker u}.
\]

In general, there just exists an isomorphism $\shF|_{X-Z} \xrightarrow{\simeq }\shF'|_{X-Z}$. But this extends to a morphism $\varphi \colon \shI \shF \to \shF'$ for some ideal $\shI \subset \shO_X$ defining some subscheme structure on $Z$ using that $X$ is noetherian \cite[I.6.9.17]{egaI2nd}. Then $\varphi$ as well as the inclusion $\shI \shF \hookrightarrow \shF$ are isomorphisms over $X-Z$. So the assertion follows by the previous case.
\end{proof}

\bigskip
One difficulty for gluing local resolutions is to guess the choice of the right syzygy sheaf $\shS$ and to solve the problem that the obstruction $o \in \Coho^2(\shHom(\shM, \shS))$ depends on both sheaves $\shM$ and $\shS$. 
We divide the construction of a locally free resolution of $\shM$ in two steps, where we can specify the first syzygy sheaf and dissolve this dependency.
For that we have to consider the local picture first.

\subsection*{Local resolutions}

Let $(R, \ideal m)$ be a noetherian local ring  and define  $X \coloneq \Spec R$,  $U \coloneq X -\lbrace \ideal m \rbrace$. By abuse of notation we identify coherent $\shO_X$-modules with their $R$-modules of global sections.

Let $M$ be an $R$-module of finite type that is locally free of constant rank on $U$. Then for every free resolution $\varphi \colon R^{\oplus n} \twoheadrightarrow M$ the first syzygy module $S \coloneq \ker \varphi$
is locally free on $U$ and satisfies $\det S|_U \simeq \det M^\vee|_U$. 
In particular, if $S$ has generically rank $1$, then $S|_U \simeq \det S|_U$ is determined by $M$, uniquely up to isomorphism.

In general, $S$ has higher rank, but we shall see that it is possible to choose a free submodule such that its quotient has rank $1$ and is locally free on $U$. 
In \cite[Theorem 2.1]{MR2041778} this was accomplished by invoking the Bourbaki Lemma \cite[p. 76]{MR0260715}.
It says that for a normal noetherian ring, every torsion-free module of rank $r$ has a free submodule  of rank $r-1$ such that its quotient is isomorphic to an ideal, hence has rank $1$.

Since we work over arbitrary noetherian rings, we need an appropriate generalization for non-reduced rings (see \cite{MR0429862} for torsion-free modules).

\begin{lem}[Modified Bourbaki Lemma]\label{LEM:bourbaki}
Let $k \in \bbN$,  $R$ be a noetherian ring and $M$ be a finitely generated $R$-module such that $M$ is free of rank $r \geq k$ at all primes of height $\leq k$.
Then there is a free submodule $F$ of $M$ of rank $r-k$, such that $M/F$ is free of rank $k$ at all primes of height $\leq k$.
\end{lem}
\begin{proof}
This is an application of basic element theory (c.f. \cite{MR811636}). Denote by $\mu(M)$ the minimal number of generators. A submodule $N \subset M$ is called \emph{$w$-fold basic} at a prime ideal $\ideal p \subset R$ if $\mu (M/N)_{\ideal p} \leq \mu( M_{\ideal p}) - w$.
A set of generators $x_1, \dots, x_s$ of $N$ is called \emph{basic up to height $k$} if $N$ is $\min(s, k -\height \ideal p + 1)$-fold basic in $M$ at each prime ideal $\ideal p\subset R$ of height less or equal to $k$.

Assume $r \geq k+1$. Let $x_1, \dots, x_s$ be a choice of generators of $M$, $s\geq r$. If $\ideal p\subset R$ is a prime ideal of height less or equal to $k$, then $w \coloneq \min(s,k -\height \ideal p + 1) \leq k+1$, thus $\mu(M_{\ideal p}) - w \geq r - (k+1)  \geq 0$, and we conclude that $x_1, \dots, x_s$ are basic up to height $k$. Hence, by \cite[Theorem 2.3]{MR811636}, there is a one element set $\lbrace y\rbrace$ which is basic up to height $k$ and induces a short exact sequence
\begin{equation}\label{LEM:bourbaki:EQ:1}
\ses{Ry}{M}{M/Ry}.
\end{equation}
Then for prime ideals $\ideal p\subset R$ with $ \height \ideal p \leq k$ the localizations
\[\ses{R_{\ideal p}y}{M_{\ideal p}}{(M/Ry)_{\ideal p}}\]
give rise to exact sequences
\[ \es{R_{\ideal p}y \otimes_R k(\ideal p)}{M_{\ideal p} \otimes_R k(\ideal p)}{(M/Ry)_{\ideal p}\otimes_R k(\ideal p)}.\]
By choice of $y$ holds $\mu(M/Ry)_{\ideal p} \leq \mu(M)_{\ideal p} - \min \lbrace 1, k -\height \ideal p +1 \rbrace = r - 1$
and we infer that $R_{\ideal p}y \otimes_R k(\ideal p)$ is nonzero. Thus, $y$ may serve as part of a basis for the free module $M_{\ideal p}$ and we conclude that $(M/Ry)_{\ideal p}$ is free of rank $r-1 \geq k$. 

By induction there is a free submodule $F' \subset M/Ry$ of rank $r-1-k$ such that $(M/Ry)/F'$ is free of rank $k$ at all primes of height $\leq k$. Pulling back the short exact sequence \eqref{LEM:bourbaki:EQ:1} along the inclusion $F' \hookrightarrow M/Ry$ gives a submodule $F \subset M$ that is an extension of $F'$ by $Ry$, thus free of rank $r-k$, and satisfies $M/F' \simeq (M/Ry)/F'$.
\end{proof}

Using the modified Bourbaki Lemma, we conclude that for a $2$-dimensional noetherian local ring, every free resolution of  an  $R$-module of finite type breaks up as follows.

\begin{prop}\label{PROP:local_resolutions}
Let $(R, \ideal m)$ be a noetherian local ring of dimension $2$ and $U = \Spec R - \lbrace \ideal m \rbrace$. Then for every finitely generated $R$-module $M$ that is locally free of rank $r \geq 1$ on $U$ there exists an exact diagram of finitely generated $R$-modules that are locally free of constant rank on $U$ such that $L|_U \simeq \det M^\vee|_U$ and $n \geq 0$:
\[
\xymatrix @R=.2in @C=.2in 
{
		& 0 \ar[d]		& 0 \ar[d] 			&		& 	\\
		& R^{\oplus n-r-1} \ar@{=}[r]\ar[d]	& R^{\oplus n-r-1} \ar[d]			&		&	\\
0 \ar[r]	& S  \ar[r]\ar[d]	& R^{n} \ar[r]^\psi\ar[d]^{\psi_2} 	& M \ar[r]\ar@{=}[d]	& 0	\\
0 \ar[r]	& L \ar[r]\ar[d]	& N \ar[r]^{\psi_1}\ar[d]		& M \ar[r]		& 0	\\
		& 0			& 0				&		&	\\
}
\]
\end{prop}
\begin{proof}
Every choice of a generating set for $M$ gives rise to a short exact sequence
\[ 
  \ses{S}{R^{\oplus n}}{M}, 
\]
such that $S$ is locally free in codimension $\leq 1$ of rank $n-r$. So from Lemma \ref{LEM:bourbaki} we obtain a short exact sequence
\[ 
  \xses{R^{\oplus n-r-1}}{}{S}{p}{L}, 
\]
such that $L$ is locally free in codimension $\leq 1$ of rank $1$. Finally, the pushout of the first exact sequence by $p$ gives the desired commutative diagram.
\end{proof}

The upshot of the previous proposition is that the surjection $\psi$ decomposes as two surjections $\psi = \psi_1 \circ \psi_2$, such that $\ker \psi_2$ is free and $\ker \psi_1$ is a coherent extension of $\det M^\vee|_U$. 
Here, the module $N$ is not free in general, but has projective dimension $\leq 1$. 

\begin{defi}
Let $X$ be a noetherian scheme. We say that a coherent $\shO_X$-module $\shN$ has property $F_k$, or is \emph{free in codimension $\leq k$} if $\shN_x$ is free for all $x \in X$ with $\dim \shO_{X,x} \leq k$ and if $\pd(\shN_x)\leq 1$ otherwise.
\end{defi}

\begin{rem}
If $X$ is a $2$-dimensional Cohen-Macaulay scheme then the Auslander-Buchsbaum formula implies that a coherent sheaf satisfies $F_k$ if and only if it has locally finite projective dimension and property $S_k$ is fulfilled.
\end{rem}

\subsection*{Global resolutions}

Building on the ideas of Proposition \ref{PROP:local_resolutions} we seek to construct a locally free resolution $\shE \twoheadrightarrow \shM$ by constructing two surjections $\psi_1 \colon \shN \twoheadrightarrow \shM$, $\psi_2 \colon \shE \twoheadrightarrow \shN$, where $\shN$ satisfies $F_1$. These surjections arise as an extension of $\shM$ by a modification of $\det \shM^\vee$ respectively by an extension of $\shN$ by a locally free sheaf.

\begin{prop}\label{PROP:resolve_coherents_by_F1sheaves}
Let $X$ be a $2$-dimensional scheme, and $\shM$ a coherent sheaf that is locally free of constant rank outside a closed set $Z \subset X$ with $\codim(Z,X)=2$, and denote by $\shF$ some chosen coherent extension of $\det \shM^\vee|_{X-Z}$.

Then there exists an obstruction
$ o \in \Coho^2(X, \shHom(\shM, \shF)),$
whose vanishing is necessary and sufficient for the existence of a short exact sequence of coherent sheaves
\[
\ses{\shL}{\shN}{\shM},
\]
where
\begin{enumerate}
 \item $\shL$ is some coherent extension of $\det \shM^\vee|_{X-Z}$, possibly different from $\shF$,
 \item $\shN$ satisfies $F_1$.
\end{enumerate}
\end{prop}
\begin{proof}
Denote by $r$ the rank of $\shM$. 
By Proposition \ref{PROP:local_resolutions} there exists for each $z\in Z$ an extension
\[
\gamma_z \colon  \quad \ses{L_z}{N_z}{M_z},
\]
such that $L_z|_{\Spec \shO_{X,z}-\lbrace z\rbrace} \simeq \det \shM^\vee|_{\Spec \shO_{X,z}-\lbrace z\rbrace}$
and $\pd(N_z) \leq 1$.
Then the family $L_z$, ${z \in Z}$, and $\det \shM^\vee|_{X-Z}$ glue to a coherent $\shO_X$-module $\shL$, i.e. $\shL_z \simeq L_z$ and $\shL|_{X-Z} \simeq \det \shM^\vee|_{X-Z}$, and the extension classes  glue to a global section $\gamma$ of $\shExt^1(\shM, \shL)$.

Since $\shHom(\shM,\shL)|_{X-Z} \simeq \shHom(\shM, \det \shM^\vee)|_{X-Z}$ we deduce from the hypothesis and Lemma \ref{LEM:coho_almost_isom} that $\Coho^2(X, \shHom(\shM,\shL)) \simeq \Coho^2(X, \shHom(\shM,\shF))$.
Thus we can identify the obstruction $o(\gamma ) \in \Coho^2(X, \shHom(\shM, \shL))$ for gluing the extensions $\gamma$ with an element $o \in \Coho^2(X, \shHom(\shM, \shF))$. If the obstruction vanishes, we obtain a coherent $\shO_X$-module $\shN$, which satisfies all asserted properties.
\end{proof}

We focus next on locally free resolutions of coherent sheaves that satisfy $F_1$. If such a resolution exists, then its first syzygy sheaf $\shS$ is locally free. Conversely, we deduce by a similar argument as in the proof of \ref{PROP:resolve_coherents_by_F1sheaves} the following existence result about \emph{locally free} resolutions.

\begin{prop}\label{PROP:resolve_F1sheaves_by_vectorbundles}
Let $X$ be a $2$-dimensional scheme and  $\shM$ a coherent sheaf satisfying $F_1$.
Then for every vector bundle $\shS$ of constant rank there exists an obstruction 
$o \in \Coho^{2}(X, \shHom(\shM, \shS^{\oplus m}))$ for some $m >0$, whose vanishing is necessary and sufficient for the existence of a locally free resolution
\[ \ses{\shS^{\oplus m}}{\shE}{\shM}. \]
\end{prop}

So up to cohomological obstructions Proposition \ref{PROP:resolve_F1sheaves_by_vectorbundles}, \ref{PROP:resolve_coherents_by_F1sheaves} and \ref{PROP:almost_ample_resolver} enable us to construct locally free resolutions of coherent sheaves on surfaces.

\begin{problem}
Proposition \ref{PROP:resolve_F1sheaves_by_vectorbundles} also holds for noetherian Deligne-Mumford stacks. However, it is not clear to the author if this is true for Proposition \ref{PROP:resolve_coherents_by_F1sheaves} since its proof depends on the fact that local sections are defined on open subsets and not just on  \'etale neighborhoods.
\end{problem}

\section{Cohomologically ample families of coherent sheaves}
\label{SEC:cample}
In this section we shall investigate families of coherent sheaves with a partial cohomological vanishing condition.
This enables us to control the cohomological constructions for gluing resolutions later on in the case of surfaces.

\begin{defi}
Let $X$ be a scheme that is proper over $A$  and let $d >0$ be an integer with $d \leq \dim(X)-1$.
A family of coherent sheaves $\shE_n$, $n \in \bbN$, is called a \emph{(cohomologically) $d$-ample family} if 
for every coherent sheaf $\shM$ there exists an $n_0 \in \bbN$ such that for all $n \geq n_0$  and $i \geq d+1$ holds
$
\Coho^i(X, \shE_n \otimes \shM) = 0.
$

A coherent sheaf $\shE$ is called \emph{cohomologically $d$-ample} if $\shE^{\otimes n}$, $n \in \bbN$, is a cohomologically $d$-ample family.
\end{defi}

If $\shE$ is a vector bundle, this cohomological vanishing condition was studied in \cite{MR0466647}, \cite{MR1618632}. For a recent treatment of line bundles we refer to \cite{Totaro:q_ample_linebundles:preprint}.

\medskip

We intend to show that the dual of an almost anti-ample sheaf is $(\dim(X)-1)$-ample (Corollary \ref{COR:dual_of_almost_antiample_is_cample}) and that $(\dim(X)-1)$-ampleness is preserved and reflected by alterations (Proposition \ref{PROP:alterations}). Recall that an \emph{alteration} is a proper morphism of schemes that is generically finite.

\begin{ex}
 An invertible sheaf on $X$ is ample if and only if it is $0$-ample.
\end{ex}

%

We will frequently exploit the fact that $d$-ampleness does not change if the sheaves of the family in question or if the scheme itself is modified at closed subsets of dimension $\leq d -1$.
The following is a direct consequence of Lemma \ref{LEM:coho_almost_isom}.

\begin{lem}\label{LEM:invariance_sheaf_mod}
Let $X$ be a scheme that is proper over $A$ and let $\shE_n$, $\shE'_n$, $n \in \bbN$ be families of coherent sheaves such that for every $n \in \bbN$, $\shE_n$ and $\shE'_n$ are isomorphic outside a closed subscheme $Z \subset X$ with $\dim(Z)\leq d - 1$.
Then $(\shE_n)$ is $d$-ample if and only if $(\shE'_n)$ is $d$-ample.
\end{lem}

Serre's Vanishing Theorem on projective birational models provides plenty of $d$-ample coherent sheaves:

\begin{prop}\label{PROP:pushforward_ample_is_dample}
Let $f \colon X' \to X$ be a morphism of proper $A$-schemes and suppose that $X'$ has an ample line bundle $\shL'$. Suppose that $f$ is an isomorphism away from a closed subset $Z \subset X$ with $\dim(Z) \leq d-1$. Then $\pf f \shL$ is $d$-ample.
\end{prop}
\begin{proof}
Let $\shM$ be a coherent sheaf on $X$ and $i \geq d +1$.
Then by Lemma \ref{LEM:coho_almost_isom} holds $\Coho^i(X,(\pf f \shL')^{\otimes n} \otimes \shM) \simeq \Coho^i(X, \pf f (\shL'^{\otimes n} \otimes \pb f \shM))$  because $\dim(Z) \leq d -1 \leq i-2$.
The latter group vanishes for $n \gg0$ using Grothendieck's spectral sequence and  Serre's Vanishing Theorem because $\shL'$ is ample, \emph{a fortiori} $f$-ample.
\end{proof}

It follows that almost anti-ample coherent sheaves (c.f. Proposition \ref{PROP:almost_ample_resolver}) satisfy a cohomological vanishing condition for the top cohomology:
\begin{cor}\label{COR:dual_of_almost_antiample_is_cample}
Let $X$ be a scheme of dimension $d \geq 1$ that is proper over $A$. Then for every almost anti-ample coherent sheaf $\shF$  the dual $\shF^\vee$ is $(d-1)$-ample.
\end{cor}
\begin{proof}
By definition there exists a proper morphism $f \colon X' \to X$ which is an isomorphism over a thick open subset $V \subset X$ and whose domain carries an ample line bundle $\shL'$ such that $\shF^\vee|_V\simeq \pf f \shL'|_V$. Since $\dim(X-V) \leq d - 2$, the assertion follows from Proposition \ref{PROP:pushforward_ample_is_dample} and Lemma \ref{LEM:invariance_sheaf_mod}.
\end{proof}

Concerning pullbacks along finite maps, we will see next that $d$-ampleness behaves as usual ampleness.
\begin{prop}\label{PROP:permanence_finite_maps}
Let $f \colon X\to Y$ be a \emph{finite} map of schemes that are proper over $A$, let $\shE_n$, $n \in \bbN$, be a family of coherent $\shO_Y$-modules and $d$ an integer with $0 \leq d < \dim(X)$.
\begin{enumerate}
 \item \label{PROP:permanence_finite_maps:IT:1} If $(\shE_n)$ is $d$-ample, then $(\pb f \shE_n)$ is $d$-ample.
 \item \label{PROP:permanence_finite_maps:IT:2} Suppose that $f$ is a surjective nilimmersion. Then $(\shE_n)$ is $d$-ample if and only if $(\pb f \shE_n)$ is  $d$-ample.
\end{enumerate}
\end{prop}
\begin{proof}
Let us prove \ref{PROP:permanence_finite_maps:IT:1} first. Given a coherent $\shO_X$-module $\shM$ the projection formula  $\pf f \shM \otimes \shE_n = \pf f (\shM \otimes \pb f \shE_n)$  holds since $\pf f$ is exact.
Therefore it induces for all $i, n \geq 0$ an isomorphism of abelian groups
$\Coho^i(X, \shM \otimes \pb f \shE_n) \simeq \Coho^i(Y, \pf f \shM \otimes \shE_n)$
which proves the assertion.
\medskip

So let us prove  \ref{PROP:permanence_finite_maps:IT:2} next. By \ref{PROP:permanence_finite_maps:IT:1} it suffices to check that the condition is sufficient. The closed immersion $f\colon X \hookrightarrow Y$ is given by a nilpotent coherent ideal $\shI \subset \shO_Y$ since $Y$ is noetherian. We may assume that $\shI^2=0$ by factoring $f$. Let $\shM$ be a given coherent $\shO_Y$-module.
Then applying $\cdot \otimes_{\shO_Y} \shM \otimes_{\shO_Y} \shE_n$ to the short exact sequence
$\ses{\shI}{\shO_Y}{\shO_X}$
we obtain a short exact sequence
\[
\ses{\shI\shM \otimes_{\shO_Y} \shE_n}{\shM \otimes_{\shO_Y} \shE_n}{\shO_X \otimes_{\shO_Y} \shM \otimes_{\shO_Y} \shE_n}
\]
which in turn gives rise to  an exact sequence
\[
 \es{\Coho^i(Y,\shI\shM \otimes_{\shO_Y} \shE_n)}{\Coho^i(Y,\shM \otimes_{\shO_Y} \shE_n)}{\Coho^i(X, \shM|_X \otimes_{\shO_X} \shE_n|_X)}.
\]
Since $\shI^2=0$ the $\shO_Y$-module $\shI \shM$ carries the structure of an $\shO_Y/\shI$-module so that we may identify the left hand side with $\Coho^i(X,\shI\shM \otimes_{\shO_X} \shE_n|_X)$. Hence, if $(\shE_n|_X)$ is $d$-ample, then also $(\shE_n)$ is $d$-ample.
\end{proof}

Finally, we show that $(\dim(X)-1)$-ampleness is preserved and reflected by pullback along alterations.

\begin{prop}\label{PROP:alterations}
Let $f \colon X \to Y$ be an alteration of schemes of dimension $d \geq 1$ that are proper over $A$. Let $\shE_n$, $n \in \bbN$, be a family of coherent $\shO_Y$-modules. Then $(\shE_n)$ is $(d-1)$-ample if and only if $(\pb f \shE_n)$ is  $(d-1)$-ample.
\end{prop}
\begin{proof}
By Proposition \ref{PROP:permanence_finite_maps}\ref{PROP:permanence_finite_maps:IT:2} we may assume that $Y$ is reduced, so $f$ has schematically dense image. Applying Stein factorization to $f$ it suffices to consider the case where $f$ is either finite or a Stein morphism.

\medskip

\emph{Case (i): $f$ finite.} 
It suffices to check that the condition is sufficient by Proposition \ref{PROP:permanence_finite_maps}.\ref{PROP:permanence_finite_maps:IT:1}.
Let $\varphi \colon \shO_Y \hookrightarrow \pf f \shO_X$ be the natural map and put $\shC \coloneq \coker \varphi$. 
Since $f$ is generically flat and finitely presented, $\shC$ is generically locally free. 
Hence $\varphi$ is generically split injective. 
Using the identification $\shHom(\pf f \shO_X, \shM) = \pf f \us f \shM$ we conclude that the transpose $\varphi^t \colon \pf f \us f \shM \to \shM$ is generically split surjective. 
It follows that $\varphi^t \otimes 1 \colon \pf f \us f \shM \otimes \shE_n \to \shM \otimes \shE_n$ is generically split surjective.
Taking cohomology gives therefore a surjection $\Coho^d(Y, \pf f \us f \shM \otimes \shE_n) \twoheadrightarrow \Coho^d(Y, \shM \otimes \shE_n)$.
The left hand side is isomorphic to $\Coho^d(Y, \pf f (\us f \shM \otimes \pb f \shE_n )) \simeq \Coho^d(X, \us f \shM \otimes \pb f \shE_n)$ using the projection formula which holds since $\pf f$ is exact. 
So if $(\pb f \shE_n)$ is $(d-1)$-ample, so too is $(\shE_n)$.
\medskip

\emph{Case (ii): $f$ Stein morphism.} Let us first verify that the condition is necessary, so let $\shM$ a given coherent $\shO_X$-module.  
Since $f$ is a proper birational map, we know that $\Coho^d(X, \shM \otimes \pb f \shE_n) \simeq \Coho^d(Y, \pf f (\shM \otimes \pb f \shE_n))$ as abelian groups by Lemma \ref{LEM:topcoho_pushforward}.
The latter is isomorphic to $\Coho^d(X, \pf f \shM \otimes \shE_n)$ by Lemma \ref{LEM:coho_almost_isom} because $f$ is an isomorphism  over all points of codimension $\leq 1$. 
This shows that if $(\shE_n)$ is $(d-1)$-ample, so too is $(\pb f \shE_n)$.

Let us next show that the condition is sufficient. By case (i) we may assume that $Y$ and $X$ are integral. 
Denote by $A$ the base ring and let $p \colon Y \to S \coloneq \Spec A$ be the structure map. 
By applying Stein factorization to $p$ and replacing $A$ we may assume that $p$ is a Stein morphism. 
For a given coherent $\shO_Y$-module $\shM$ we have to show that the $A$-module $\Coho^d(Y, \shM \otimes \shE_n)=\Coho^0(S, R^d \pf p (\shM \otimes \shE_n))$ vanishes for $n \gg0$.
The coherent sheaf $R^d \pf p (\shM \otimes \shE_n)$ is zero if its stalk at all closed points vanishes. 
So we may assume that $\dim p^{-1}(s) \geq d$ for some closed point $s \in S$ by the Theorem on Formal Functions. 
But then $\dim p^{-1}(s)=d$ and since $Y$ is irreducible it follows that for the generic point  $\eta \in Y$ holds $p(\eta)=s$.
Using that $s$ is is closed and $p$ is proper we infer $p(Y)=s$. So $S$ consists of a single point. Since $Y$ is reduced and $p$ is Stein it follows that $S$ is the spectrum of a field. Thus $X$ and $Y$ are algebraic schemes.

It follows that the normalizations of $X' \to X$ and $Y' \to Y$ are \emph{finite} maps. Therefore we may assume that $X$ and $Y$ are normal, proper over a field and that $f$ is birational by case (i). 
Thus, the assertion follows from the succeeding lemma.
\end{proof}
\begin{lem}
Let $f \colon Y \to X$ be a birational map of $d$-dimensional, normal schemes that are proper over a field. Then for every coherent $\shO_X$-module $\shF$ there exists a surjection 
\[\Coho^d(Y, \shHom(\pb f \omega_X, \omega_Y) \otimes \pb f \shF ) \twoheadrightarrow \Coho^d(X, \shF).\]
\end{lem}
\begin{proof}
By Serre duality there is an isomorphism $\Coho^d(X, \shF)^\vee \simeq \Hom(\shF, \omega_X)$.
Since $X$ is normal the dualizing module $\omega_X$ is torsion-free and generically invertible and the birationality of $f$ implies the injectivity of the natural map
\(
\Hom( \shF, \omega_X) \rightarrow 
\Hom(\pb f \shF,  \pb f \omega_X) 
\rightarrow \Hom(\pb f \shF, \shHom( \shHom( \pb f \omega_X,\omega_Y),\omega_Y)).
\)
By adjunction the group on the right is isomorphic to
\(
\Hom( \pb f \shF \otimes \shHom( \pb f \omega_X,\omega_Y), \omega_Y)
\)
which is isomorphic to
\(
\Coho^d(Y, \pb f \shF \otimes \shHom( \pb f \omega_X,\omega_Y))^\vee
\)
using Serre duality.
\end{proof}

\begin{lem}\label{LEM:topcoho_pushforward}
Let $f \colon X \to Y$ be a proper birational Stein morphism of noetherian schemes of dimension $d$. Then $\Coho^d(X, \shF) \simeq \Coho^d(Y, \pf f \shF)$ for every  coherent $\shO_X$-module $\shF$.
\end{lem}
\begin{proof}
Let $y\in Y$ with $\dim \shO_{y} \geq 1$. Then $\dim f^{-1}(x) \leq \dim \shO_x-1$ since $f$ is birational and of finite type.
Hence $(R^q \pf f \shM)_y = 0$ for all $q \geq \dim \shO_x$ by the Theorem on Formal Functions.
It follows that $\codim( \Supp R^q \pf f \shM, Y) \geq q+1$ and this implies $\Coho^p(Y, R^q \pf f \shM) =0$ for all $p\geq 0$, $q \geq 1$ with $p+q \geq d$. Applying the Grothendieck spectral sequence settles the result.
\end{proof}

\section{Existence of positive vector bundles on non-projective surfaces}
\label{SEC:existence_of_cample_families_of_vectorbundles}

In this section we construct a $1$-ample family of vector bundles $\shE_n$, $n \in \bbN$, of rank $2$ on an arbitrary surface $X$ that is proper over a $A$ (see Theorem \ref{THM:existence_of_1ample_rank2_family}). 

\medskip

Let us first discuss descent conditions of vector bundles for a proper birational map of surfaces.
For that we have to introduce some terminology. Let $f \colon X \to Y$ be a proper birational morphism. 
The closed subscheme $B \subset Y$ given by the conductor ideal $\Ann_{\shO_Y}\coker(\shO_Y \to \pf f \shO_X)$ is called the \emph{branching subscheme}. The union of all integral $1$-dimensional closed subschemes contracted by $f$ is called the \emph{exceptional curve} $E \subset X$.

\begin{lem}\label{LEM:vector_bundle_descent}
Let $f \colon X \to Y$ a birational morphism of $2$-dimensional noetherian schemes that satisfies the following conditions:
\begin{enumerate}
 \item \label{LEM:vector_bundle_descent:IT:1} The branching subscheme $B \subset Y$ is $0$-dimensional.
 \item \label{LEM:vector_bundle_descent:IT:2} There exists an effective Cartier divisor $D \subset X$ that is supported on the exceptional curve $E \subset X$ such that $\shO_E(-D)$ is ample.
\end{enumerate}
Then for every $r \in \bbN$ there is an $m \in \bbN$ such that for every vector bundle $\shE$ on $X$ of rank $r$ whose restriction $\shE|_{mE}$ is trivial, there exists a vector bundle $\shF$ on $Y$ of rank $r$ with $\pb f \shF \simeq \shE$.
\end{lem}
\begin{proof}
Let $X \overset{f_0} \to Y_0 \overset{\iota}\hookrightarrow Y$ the the factorization of $f$ over its schematic image.
Then every vector bundle on $Y_0$ lifts to a vector bundle on $Y$.
 To see this we may assume that the nilimmersion $Y_0 \hookrightarrow Y$ is given by a coherent ideal $\shI \subset \shO_Y $ with $\shI^2=0$, so that $\shI$ carries the structure of an $\shO_{Y_0}$-module.
Then the obstruction of lifting a locally free $\shO_{Y_0}$-module $\shF_0$ is an element of $\Coho^2(Y_0, \shI \otimes_{\shO_{Y_0}} \shEnd(\shF_0))$ \cite[Theorem 5.3]{MR2223409}. But this group vanishes as the support of $\shI$ has dimension $\leq 1$.
The branching subscheme $B_0 \subset Y_0$ of $f_0$ equals $\iota^{-1}(B)$ and the $f$- and $f_0$-exceptional curve coincides, so the conditions \ref{LEM:vector_bundle_descent:IT:1}, \ref{LEM:vector_bundle_descent:IT:2} on $f$ carry over to $f_0$.
So we may assume that $f$ has schematically dense image. Then the same proof as for \cite[1.2]{MR2041778} applies.
\end{proof}

First, we shall prove in Proposition \ref{PROP:existence_vectorbundles_on_projective_surface} a stronger variant of Theorem \ref{THM:existence_of_1ample_rank2_family} in the projective case. 

\begin{lem}\label{LEM:trivial_bundle_representation_on_curve}
Let $\shL$ be an ample line bundle on a $1$-dimensional scheme $X$. Then there exists a short exact sequence
\[ \ses{\shO_X^{\oplus 2}}{(\shL^n)^{\oplus 2}}{\shO_D} \]
for some effective Cartier divisors $D \subset X$ with $\shO_X(D) \simeq \shL^{2n}$ and $n >0$.
\end{lem}
\begin{proof}
Let $\beta \in \Coho^0(X, \shL^b)$, $b \gg 0$, be a global section that is non-zero at all associated points of $X$ \cite[4.5.4]{egaII}. Then $\beta \colon \shO_X \to \shL^b$ is injective and $\coker \beta \simeq \shO_B$ for the  Cartier divisor $B \coloneq V(\beta)$ which satisfies $\shO_X(B) \simeq \shL^b$. 
Let $\alpha \in \Coho^0(X, \shL^a)$, $a=bc$, $c \gg0$, a second global section that is nonzero over  $B \cup \Ass(X)$. Then $\alpha \colon \shO_X \to \shL^a$ is also injective, $\coker t \simeq \shO_A$ for the  Cartier divisor $A \coloneq V(\alpha) \subset X$ that is  disjoint from $B$ and  satisfies $\shO_X(A)$. 
Then $\alpha^{\otimes c} \oplus \beta \colon \shO_X^{\oplus 2} \to (\shL^a)^{\oplus 2}$ is also injective and its cokernel is isomorphic to $\shO_{cB+A}$.
\end{proof}

\begin{lem}\label{LEM:interpolate}
Let $X$ be a $1$-dimensional scheme that is proper over $A$ and has an ample line bundle $\shL$.  Then for every $a_0, b_0 \in \Coho^0(Z,\shO_Z)$ where $Z \subset X$ is a discrete closed subscheme disjoint from $\Ass(X)$, there exist $a,b \in \Coho^0(X,\shL^n)$ for some $n >0$ and a regular section $s \in \Coho^0(Z, \shL^n|_Z)$ such that
\begin{enumerate}
 \item \label{LEM:interpolate:IT:1} $a_0 s = a|_Z$ and  $b_0 s  = b|_Z$,
 \item \label{LEM:interpolate:IT:2} $X = Z \cup X_a \cup X_b$.
\end{enumerate}
\end{lem}
\begin{proof}
By enlarging $Z$ and extending $a_0$, $b_0$ by $1$ we may assume that $a_0$ and $b_0$ are non-zero at all embedded points and on each irreducible component of $X$.
Using that $\shL$ is ample and replacing $\shL$ with a multiple, there exists a section $s \in \Coho^0(X,\shL)$  such that $X_s \subset X$ is an affine open neighborhood of $Z$. 
We can choose a section $a' \in \Coho^0(X_s,\shO_{X_s})$ that satisfies $a'|_Z = a_0$.
Then $a' s^{\otimes k}$ extends to a section $a'' \in \Coho^0(X, \shL^{k})$ for some $n >0$ \cite[6.8.1]{egaI2nd}. 
It follows that $a|_Z = a_0 s^{\otimes k}|_Z$ and $V(a) \subseteq V(a') \cup V(s) \subset X$ is discrete.

Using that $\shL$ is ample and $X$ is proper over $A$, for $Y \coloneq V(a) - Z$ the exact sequence
\[
 \es{ \Coho^0(X, \shL^{l})}
  {\Coho^0(Y \cup Z, \shL^{l}|_{ Y \cup Z} )}
  {\Coho^1(X, \shI_{Y \cup Z} \otimes \shL^{l})}
\]
implies that for some $l>n$ there is a section $t \in \Coho^0(X,\shL^{l})$ which is non-zero near $Y$ and satisfies $t|_Z = s^{\otimes l}|_Z$. 
It follows that $X_t \subset X$ is an open neighborhood of $Y\cup Z$, hence dense and affine.
Now choose a section $b' \in \Coho^0(X_t, \shO_{X_t})$ that is non-zero near $Y$ and satisfies $b'|_Z = b_0$. Then $b' t^{\otimes m}$ extends to a section $b \in \Coho^0(X, \shL^{lm})$ for some $m>0$, satisfying \( b|_Z = b_0  t^{\otimes m}|_Z = b_0  s^{\otimes lm}|_Z\) and $V(b) \subset  X_a \cup Z$.

With $d \coloneq ml-k>0$ follows
\(  a \otimes s^{\otimes d}|_Z = a_0  s^{\otimes k}|_Z \otimes  s^{\otimes ml-k}|_Z = a_0  s^{ \otimes lm}|_Z\) which implies \ref{LEM:vector_bundle_descent:IT:1} for $n=klm$ and $s$ replaced by $s^{\otimes lm}$.
Finally \ref{LEM:vector_bundle_descent:IT:2} also holds because $X_{a \otimes s^{\otimes d}} \cup X_b = (X_a \cap X_s) \cup X_b = X_a \cup X_b \supset X-Z$.
\end{proof}

\begin{prop}\label{PROP:existence_vectorbundles_on_projective_surface}
Let $X$ be a surface that is proper over $A$, satisfies $S_1$ and has an ample line bundle $\shL$. 
Then for every Weil divisor $Y \subset X$ there exists a vector bundle $\shE$ on $X$ whose restriction $\shE|_Y$ is trivial and that fits in a short exact sequence
\begin{equation}\label{PROP:existence_vectorbundles_on_projective_surface:EQ:0}
 \ses{\shE}
  {(\shL^{n})^{\oplus 2}}
  {\shL^m|_X},
\end{equation}
where $C \subset X$ is an effective Cartier divisor with $\shO_X(C) \simeq \shL^{2n}$ and $m > n> 0$.
\end{prop}
\begin{proof}
By enlarging $Y$ we may assume that every irreducible component of $X$ contains some component of $Y$.
Let $\shI_Y \subset \shO_X$ be the coherent ideal defining $Y \subset X$ and suppose that $\Coho^1(X, \shI_Y \otimes \shL^a)=0$ for all $a \geq 1$ by replacing $\shL$ with a suitable multiple.

Replacing $\shL$ again with a suitable multiple, there is a short exact sequence by Lemma \ref{LEM:trivial_bundle_representation_on_curve}:
\begin{equation}\label{PROP:existence_vectorbundles_on_projective_surface:EQ:1}
 \ses{\shO_Y^{\oplus 2}}
  {\shL|_Y^{\oplus 2}}
  {\shO_Z}
\end{equation}
for some Cartier divisor $Z \subset Y$ satisfying $\shO_Y(Z) \simeq \shL^2|_Y$. 
Applying $\cdot \otimes \shL^{-1}|_Y$ induces a short exact sequence 
\begin{equation}\label{PROP:existence_vectorbundles_on_projective_surface:EQ:2}
 \xses{(\shL^{-1})|_Y^{\oplus 2}}{}{\shO_Y^{\oplus 2}}{\varphi}{\shO_Z}.
\end{equation}

\begin{claim*}
There exists an effective Cartier divisor $C \subset X$ with $C \cap Y = Z$ and $\shO_X(C)= \shL^2$.
\end{claim*}
By assumption on $\shL$ the right-hand side of the exact sequence
\[
 \es{\Coho^0(X, \shL^2)}
  {\Coho^0(Y, \shL^2|_Y)}
  {\Coho^1(X, \shI_Y \otimes \shL^2)}
\]
vanishes, so that the section $z' \in \Coho^0(Y, \shL^2|_Y)$ defining $Z \subset Y$ lifts to global section $z \in \Coho^0(X, \shL^2)$. 
Then $z$ is non-zero at all generic points of $X$ because $X_z \cap Y = Y_{z'} = Y-Z$ meets every irreducible component of $X$ by assumption on $Y$.
Since $X$ has no embedded points it follows that $z$ is a \emph{regular} section, so that $C = V(z)$ is a Cartier divisor satisfying the asserted properties.

\medskip

The map $\varphi$ is given as $(y, y') \mapsto y|_Z a_0 + y'|_Z b_0$ for two sections $a_0,b_0 \in \Coho^0(Z, \shO_Z)$ and all sections $y,y'$ of $\shO_Y$.
Using Lemma \ref{LEM:interpolate} there exists a trivialization $\sigma \colon \shO_Z \xrightarrow{\sim} \shL^n|_Z, \, 1 \mapsto s$, and two sections $a,b \in \Coho^0(C, \shL^n|_C)$ such that $a|_Z=a_0 s$, $b|_Z=b_0s$ and $X_a \cup X_a \supset X-Z$. 

Define $\Phi \colon \shO_X^{\oplus 2} \to \shL^n|_C$ by $(x,x') \mapsto x|_C a + x'|_C b$. 
Then $\Phi|_Y = \sigma \circ \varphi$. In particular, $\Phi|_Z$ is surjective. But $\Phi|_{X-Z}$ is also surjective because $X_a \cup X_a \supset X-Z$. The upshot is that $\Phi$ is a surjective map that extends $\varphi$ up to isomorphism of the codomain. 
Since $\pd \shL^n|_C \leq 1$ it follows that $\shF \coloneq \ker \Phi$ is a locally free $\shO_X$-module of rank $2$.

\begin{claim*}
For the restriction holds $\shF|_Y \simeq \shL^{-1}|_Y$.
\end{claim*}
In fact, restricting the short exact sequence
\begin{equation*}
 \xses{\shF} {\psi} {\shO_X^{\oplus 2}} {\Phi} { \shL^{n}|_C}
\end{equation*}
to $Y$ gives an exact sequence
\begin{equation*}
\xses{\shF|_Y} {\psi|_Y} {\shO_Y^{\oplus 2}} {\sigma \circ \varphi} { \shL^{n}|_Z}.
\end{equation*}
Here $\psi|_Y$ is injective because $\shF|_Y$ has no sections supported on $Z$. Therefore $\shF|_Y \simeq \shL^{-1}|_Y$ using \eqref{PROP:existence_vectorbundles_on_projective_surface:EQ:2} and $\shE \coloneq \shF \otimes \shL$ is the desired vector bundle.
\end{proof}

\begin{thm}\label{THM:existence_of_1ample_rank2_family}
 Let $X$ be a $2$-dimensional scheme that is proper over $A$. Then there exists a cohomologically $1$-ample family of vector bundles $\shE_i$, $i \in \bbN$, of rank $2$.
\end{thm}
\begin{proof}
Let us first reduce to the case that $X$ is projective over $A$.
By Lemma \ref{LEM:vector_bundle_descent} every vector bundle descends along the $S_1$-ization $X' \to X$. So we may assume that $X$ has no embedded points in light of Proposition \ref{PROP:alterations}.\ref{PROP:permanence_finite_maps:IT:2}.

By Theorem \ref{THM:thick_quasiprojectives} there exists an open subset $V \subset X$ which is quasiprojective over $A$ and $V-X$ consists of finitely many points of codimension $2$. 
Then by Nagata there exists a $V$-admissible blow up $f \colon Y \to X$ such that $Y$ is projective over $A$ \cite[2.6]{MR2356346}. Also, $Y$ has no embedded points. 

Let us verify the descent conditions of Lemma \ref{LEM:vector_bundle_descent} vor $f$.
Clearly, $f$ has $0$-dimensional branching subscheme since $V$ is thick.
Using that $Z$ is $0$-dimensional the inverse image $D \coloneq f^{-1}(Z)$ is supported on the $f$-exceptional curve $E \subset Y$. By construction of the blow-up the closed subscheme $D \subset Y$ is an effective Cartier divisor given by the $f$-ample invertible inverse image ideal $\shI_Z \cdot \shO_Y = \shO_Y(-D)$. Thus $\shO_E(D)$ is ample.
Then by Lemma \ref{LEM:vector_bundle_descent} and Proposition \ref{PROP:alterations} there exists an $m \in \bbN$ such that every $1$-ample family of vector bundles $\shF_n$, $n \in \bbN$, of rank $2$ on $Y$ descends to a $1$-ample family vector bundles $\shE_n$ of rank $2$ on $X$, if each restriction $\shF_n|_{mE}$ is trivial.

So we may assume that $X$ is projective over $A$ if we additionally show that the family $\shE_n$, $n \in \bbN$, is trivial on a given Weil divisor $Y \subset X$. 
Choose an ample $\shO_X$-module $\shL$. 
By Proposition \ref{PROP:existence_vectorbundles_on_projective_surface} there exists for every $n \in \bbN$ a vector bundle $\shE_n$ whose restriction $\shE_n|_Y$ is trivial and that fits in a short exact sequence
\[
 \ses{\shE_n}{(\shL^{a_n })^{\oplus 2}}{\shL^{b_n}|_C}
\]
for some Cartier divisor $C \subset X$ and integers $b_n > a_n > n$.
Given a coherent $\shO_X$-module $\shM$, applying $ \cdot \otimes \shM$ results in an exact sequence
\[
 0 \rightarrow 
 \shN_n \rightarrow 
 \shE_n \otimes \shM \xrightarrow{\varphi_n} 
 (\shL^{a_n })^{\oplus 2}  \otimes \shM  \rightarrow 
 \shL^{b_n}|_C  \otimes \shM 
 \rightarrow 0
\]
for some coherent $\shO_X$-module $\shN_i$ supported on $C$. Taking the associates long exact cohomology sequences gives two exact sequences
\[
\es{\Coho^2(Y, \shN_n)}
  {\Coho^2(Y,\shE_n \otimes \shM)}
  {\Coho^2(Y, \im \varphi_n)}
\]
\[
 \xes{\Coho^1(Y, \shL^{b_n}|_C  \otimes \shM)}{}
     {\Coho^2(Y, \im \varphi_n)}{\partial}{\Coho^2(Y, (\shL^{a_n })^{\oplus 2})}
\]
Using that $\shL$ is ample, we conclude that $\Coho^2(Y, \im \varphi_n)$ vanishes for $n \gg 0$.
But $\Coho^2(Y, \shN_n)$ is also zero because the support of $\shN_n$ has dimension $\leq 1$.
This shows that $\shE_n$, $n \in \bbN$, is a cohomologically $1$-ample family of vector bundles of rank $2$.
\end{proof}

\section{Proof of the resolution property for surfaces}
\label{SEC:resprop_surfaces}

Let us finally collect the preceding results to prove the resolution property for proper surfaces.
\begin{thm}\label{THM:main_theorem}
Every $2$-dimensional scheme that is proper over $A$ satisfies the resolution property.
\end{thm}
\begin{proof}
Let $\shM$ be a coherent sheaf and $x \in X$ an arbitrary point. 
By Proposition \ref{PROP:almost_ample_resolver} there exists a coherent sheaf $\shF$ that is almost anti-ample near $x$. In particular, for every $m \gg0$ there exists a map $(\shF^{\otimes m})^{\oplus n} \to \shM$ for some $n \in \bbN$ that is surjective near $x$. Therefore it suffices to find a locally free resolution for every $\shF^{\otimes m}$ and every $m \gg 0$ because $X$ is quasicompact.

By definition, $\shF$ is invertible outside finitely many points of codimension $2$.
So by Proposition \ref{PROP:resolve_coherents_by_F1sheaves} for every $m \gg 0$ there exists a surjection $\shG \to \shF^{\otimes m}$, such that $\shG$ satisfies $F_1$ because $\Coho^2(X, \shHom(\shF^{\otimes m}, (\shF^\vee)^{\otimes m})) \simeq \Coho^2(X, (\shF^\vee)^{\otimes 2m})$ (Lemma \ref{LEM:coho_almost_isom}) vanishes using that $\shF^\vee$ is  $1$-ample by Corollary \ref{COR:dual_of_almost_antiample_is_cample}.

By Theorem \ref{THM:existence_of_1ample_rank2_family} there exists a $1$-ample family $(\shE_n)$, $n \in \bbN$, of vector bundles of rank $2$. 
This implies  $\Coho^2(X, \shHom(\shG, \shE_n)) \simeq \Coho^2(X, \shG^\vee \otimes \shE_n) = 0$ for $n \gg0$.
Consequently $\shG$ admits a surjection $\shH \twoheadrightarrow \shG$ by a vector bundle $\shH$ in light of Proposition \ref{PROP:resolve_F1sheaves_by_vectorbundles}.
\end{proof}

Since the resolution property descends along immersions, it holds for all $2$-dimensional schemes $X$ that are embeddable into $2$-dimensional schemes which are proper over a noetherian base ring. 

\begin{cor}\label{COR:algebraic_surfaces}
Suppose that $A$ is a noetherian, universally catenary Jacobson ring such that each irreducible component of $\Spec A$ is equicodimensional (for example, if $A=\bbZ$ or if $A$ is a field).
Then every $2$-dimensional scheme that is separated and of finite type over $A$ satisfies the resolution property.
\end{cor}
\begin{proof}
Since $X$ is separated and of finite type over a noetherian ring, there exists a proper $A$-scheme $\overline X$ together with an open immersion $X \hookrightarrow \overline X$ which identifies $X$ as a dense open subscheme of $\overline X$ by Nagata's embedding Theorem \cite[4.1]{MR2356346}.
The assumptions on $A$ guarantee that $\dim(\overline X) = \dim (X)=2$ since $X$ is of finite type over $A$ \cite[10.6.2]{egaIV_3}.
So Theorem \ref{THM:main_theorem} implies that $\overline X$ and hence $X$ has the resolution property.
\end{proof}

\begin{rem}
In the first section we extended the resolution property from a given dense affine open subset $U_0 \subset X$ to a thick open subset $U_0 \subset U_1 \subset X$, i.e. we added all points of codimension $1$. The results of Section \ref{SEC:glueing_resolutions} can be adapted to formulate conditions to extend the resolution property to an open subset $U_1 \subset U_2 \subset X$ that contains all points of codimension $2$ and the cohomological obstructions lie then in $2^{nd}$ cohomology groups of coherent $\shO_{U_2}$-modules. With a similar argument as in the end of the proof of Theorem \ref{THM:thick_quasiprojectives} one can arrange that $U_2$ has cohomological dimension $\leq 2$. 
We believe that $U_2$ then already satisfies the resolution property up to a closed subset $Z \in X-U_1$ of codimension $\geq 3$.
For technical reasons we assumed that $X$ is $2$-dimensional and proper over a noetherian ring, so that $U_2=X$ is proper which allows to control the cohomological obstructions.
\end{rem}

\bibliographystyle{pgalpha}
\bibliography{pgbiblio}

\begin{thebibliography}{EGA I\textsubscript{2nd}}

\bibitem[Art71]{MR0289501}
M.~Artin.
\newblock On the joins of {H}ensel rings.
\newblock {\em Advances in Math.}, 7:282--296 (1971), 1971.

\bibitem[Bor67]{MR0219545}
Mario Borelli.
\newblock Some results on ampleness and divisorial schemes.
\newblock {\em Pacific J. Math.}, 23:217--227, 1967.

\bibitem[Bou65]{MR0260715}
N.~Bourbaki.
\newblock {\em \'{E}l\'ements de math\'ematique. {F}asc. {XXXI}. {A}lg\`ebre
  commutative. {C}hapitre 7: {D}iviseurs}.
\newblock Actualit\'es Scientifiques et Industrielles, No. 1314. Hermann,
  Paris, 1965.

\bibitem[BS03]{MR1970862}
Holger Brenner and Stefan Schr{\"o}er.
\newblock Ample families, multihomogeneous spectra, and algebraization of
  formal schemes.
\newblock {\em Pacific J. Math.}, 208(2):209--230, 2003.

\bibitem[BV75]{MR0429862}
Winfried Bruns and Udo Vetter.
\newblock Die {V}erallgemeinerung eines {S}atzes von {B}ourbaki und einige
  {A}nwendungen.
\newblock {\em Manuscripta Math.}, 17(4):317--325, 1975.

\bibitem[Con07]{MR2356346}
Brian Conrad.
\newblock Deligne's notes on {N}agata compactifications.
\newblock {\em J. Ramanujan Math. Soc.}, 22(3):205--257, 2007.

\bibitem[EG85]{MR811636}
E.~Graham Evans and Phillip Griffith.
\newblock {\em Syzygies}, volume 106 of {\em London Mathematical Society
  Lecture Note Series}.
\newblock Cambridge University Press, Cambridge, 1985.

\bibitem[EGA I\textsubscript{2nd}]{egaI2nd}
A.~Grothendieck and Jean~A. Dieudonn\'e.
\newblock {\em \'{E}l\'ements de g\'eom\'etrie alg\'ebrique. I.}
\newblock {D}ie Grundlehren der mathematischen {W}issenschaften. 166.
  {B}erlin-{H}eidelberg-{N}ew {Y}ork {S}pringer-{V}erlag. IX, 466 p., 1971.

\bibitem[EGA II]{egaII}
A.~Grothendieck.
\newblock \'{E}l\'ements de g\'eom\'etrie alg\'ebrique. {II}. \'{E}tude globale
  \'el\'ementaire de quelques classes de morphismes.
\newblock {\em Inst. Hautes \'Etudes Sci. Publ. Math.}, (8):222, 1961.

\bibitem[EGA III\textsubscript{1}]{egaIII_1}
A.~Grothendieck.
\newblock \'{E}l\'ements de g\'eom\'etrie alg\'ebrique. {III}. \'{E}tude
  cohomologique des faisceaux coh\'erents. {I}.
\newblock {\em Inst. Hautes \'Etudes Sci. Publ. Math.}, (11):167, 1961.

\bibitem[EGA IV\textsubscript{2}]{egaIV_2}
A.~Grothendieck.
\newblock \'{E}l\'ements de g\'eom\'etrie alg\'ebrique. {IV}. \'{E}tude locale
  des sch\'emas et des morphismes de sch\'emas. {II}.
\newblock {\em Inst. Hautes \'Etudes Sci. Publ. Math.}, (24):231, 1965.

\bibitem[EGA IV\textsubscript{3}]{egaIV_3}
A.~Grothendieck.
\newblock \'{E}l\'ements de g\'eom\'etrie alg\'ebrique. {IV}. \'{E}tude locale
  des sch\'emas et des morphismes de sch\'emas. {III}.
\newblock {\em Inst. Hautes \'Etudes Sci. Publ. Math.}, (28):255, 1966.

\bibitem[Fer03]{MR2044495}
Daniel Ferrand.
\newblock Conducteur, descente et pincement.
\newblock {\em Bull. Soc. Math. France}, 131(4):553--585, 2003.

\bibitem[Ill05]{MR2223409}
Luc Illusie.
\newblock Grothendieck's existence theorem in formal geometry.
\newblock In {\em Fundamental algebraic geometry}, volume 123 of {\em Math.
  Surveys Monogr.}, pages 179--233. Amer. Math. Soc., Providence, RI, 2005.
\newblock With a letter (in French) of Jean-Pierre Serre.

\bibitem[Jel05]{MR2146196}
Zbigniew Jelonek.
\newblock On the projectivity of threefolds.
\newblock {\em Proc. Amer. Math. Soc.}, 133(9):2539--2542 (electronic), 2005.

\bibitem[Kle66]{MR0206009}
Steven~L. Kleiman.
\newblock Toward a numerical theory of ampleness.
\newblock {\em Ann. of Math. (2)}, 84:293--344, 1966.

\bibitem[Mum70]{MR0282985}
David Mumford.
\newblock {\em Abelian varieties}.
\newblock Tata Institute of Fundamental Research Studies in Mathematics, No. 5.
  Published for the Tata Institute of Fundamental Research, Bombay, 1970.

\bibitem[Pay09]{MR2448277}
Sam Payne.
\newblock Toric vector bundles, branched covers of fans, and the resolution
  property.
\newblock {\em J. Algebraic Geom.}, 18(1):1--36, 2009.

\bibitem[Ray70]{MR0260758}
Michel Raynaud.
\newblock {\em Faisceaux amples sur les sch\'emas en groupes et les espaces
  homog\`enes}.
\newblock Lecture Notes in Mathematics, Vol. 119. Springer-Verlag, Berlin,
  1970.

\bibitem[RV04]{MR2096147}
Mike Roth and Ravi Vakil.
\newblock The affine stratification number and the moduli space of curves.
\newblock In {\em Algebraic structures and moduli spaces}, volume~38 of {\em
  CRM Proc. Lecture Notes}, pages 213--227. Amer. Math. Soc., Providence, RI,
  2004.

\bibitem[Sch82]{MR676049}
Hans-Werner Schuster.
\newblock Locally free resolutions of coherent sheaves on surfaces.
\newblock {\em J. Reine Angew. Math.}, 337:159--165, 1982.

\bibitem[Sch99]{MR1726231}
Stefan Schr{\"o}er.
\newblock On non-projective normal surfaces.
\newblock {\em Manuscripta Math.}, 100(3):317--321, 1999.

\bibitem[SGA 6]{sga6}
{\em Th\'eorie des intersections et th\'eor\`eme de {R}iemann-{R}och}.
\newblock Lecture Notes in Mathematics, Vol. 225. Springer-Verlag, Berlin,
  1971.
\newblock S{\'e}minaire de G{\'e}om{\'e}trie Alg{\'e}brique du Bois-Marie
  1966--1967 (SGA 6), Dirig{\'e} par P. Berthelot, A. Grothendieck et L.
  Illusie. Avec la collaboration de D. Ferrand, J. P. Jouanolou, O. Jussila, S.
  Kleiman, M. Raynaud et J. P. Serre.

\bibitem[Som78]{MR0466647}
Andrew~John Sommese.
\newblock Submanifolds of {A}belian varieties.
\newblock {\em Math. Ann.}, 233(3):229--256, 1978.

\bibitem[Ste98]{MR1618632}
Frauke Steffen.
\newblock A generalized principal ideal theorem with an application to
  {B}rill-{N}oether theory.
\newblock {\em Invent. Math.}, 132(1):73--89, 1998.

\bibitem[SV04]{MR2041778}
Stefan Schr{\"o}er and Gabriele Vezzosi.
\newblock Existence of vector bundles and global resolutions for singular
  surfaces.
\newblock {\em Compos. Math.}, 140(3):717--728, 2004.

\bibitem[Tot04]{MR2108211}
Burt Totaro.
\newblock The resolution property for schemes and stacks.
\newblock {\em J. Reine Angew. Math.}, 577:1--22, 2004.

\bibitem[Tot10]{Totaro:q_ample_linebundles:preprint}
Burt Totaro.
\newblock Line bundles with partially vanishing cohomology.
\newblock {\em arXiv:1007.3955v1 [math.AG]}, 2010.

\bibitem[Voi02]{MR1902630}
Claire Voisin.
\newblock A counterexample to the {H}odge conjecture extended to {K}\"ahler
  varieties.
\newblock {\em Int. Math. Res. Not.}, (20):1057--1075, 2002.

\bibitem[W{\l}o93]{MR1227474}
Jaros{\l}aw W{\l}odarczyk.
\newblock Embeddings in toric varieties and prevarieties.
\newblock {\em J. Algebraic Geom.}, 2(4):705--726, 1993.

\bibitem[W{\l}o99]{MR1683254}
Jaros{\l}aw W{\l}odarczyk.
\newblock Maximal quasiprojective subsets and the {K}leiman-{C}hevalley
  quasiprojectivity criterion.
\newblock {\em J. Math. Sci. Univ. Tokyo}, 6(1):41--47, 1999.

\end{thebibliography}
\end{document}